\documentclass[a4paper,11pt]{article}
\usepackage{amsfonts,amsmath,amssymb,amsthm,enumerate,bm,dsfont,a4,}
\usepackage[a4paper,text={6.25in,9.0in},centering]{geometry}
\allowdisplaybreaks[4] 
\theoremstyle{plain}
\newtheorem{theorem}{\noindent\bf Theorem}[section]
\newtheorem{corollary}[theorem]{\noindent\bf Corollary}

\newtheorem{lemma}[theorem]{\noindent\bf Lemma}

\theoremstyle{remark}

\newtheorem{remark}[theorem]{\noindent\bf Remark}

\numberwithin{equation}{section}

\def\be{\begin{eqnarray}}%%
\def\ee{\end{eqnarray}}%%
\def\ben{\begin{eqnarray*}}%%
\def\een{\end{eqnarray*}}%%
\def\benum{\begin{enumerate}}%%
\def\eenum{\end{enumerate}}%%
\newcommand{\lr}{\left(}
\newcommand{\rr}{\right)}
\newcommand{\lek}{\left[}
\newcommand{\rek}{\right]}
\newcommand{\lge}{\left\{ }
\newcommand{\rge}{\right\} }

\newcommand{\vertiii}[1]{{\left\vert\kern-0.25ex\left\vert\kern-0.25             ex\left\vert #1 \right\vert\kern-0.25ex\right\vert     
            \kern-0.25ex\right\vert}}
\title{\bf Dunkl-Hausdorff operator on weighted Orlicz spaces}
\author{Megha Madan, Arun Pal Singh,  Monika Singh \footnote{the corresponding author}}

\date{}

\begin{document}
\maketitle

\begin{abstract}
	 We characterize the modular and norm inequalities for the  Dunkl-Hausdorff operator defined on non-negative non-increasing functions in the framework of the weighted Orlicz spaces.	
\end{abstract}
	
	\vspace{15pt}
	

   \section{Introduction and Preliminaries}
   
Dunkl theory generalizes classical Fourier analysis on
   $\mathbb{R}^n.$ It started with the Dunkl’s seminal work in \cite{dk} and was further developed by several mathematicians. Liflyand and M\'oricz in \cite{lm1, lm2} studied respectively the boundedness of Hausdorff operator and its commuting relation with Hilbert transform on real Hardy spaces. In the frame of extending the results of Liflyand and M\'oricz to the context of Dunkl  theory, Daher and Saadi in \cite{dsa} introduced and studied  the  Dunkl-Hausdorff operator on weighted Sobolev spaces $L^p_k(\mathbb{R})$ and Hardy spaces $\mathcal{H}^1_k(\mathbb{R}).$ The boundedness of the Dunkl-Hausdorff operator on Lebesgue spaces in one dimension as well as two dimensions was studied in \cite{jfj}. Also, the Dunkl-Hausdorff operator has been studied for monotone functions and monotone weights on Lebesgue spaces in \cite{jj} and \cite{jjsm}. In this paper we study Dunkl-Hausdorff operator on Orlicz spaces.  
   
   Orlicz spaces were introduced by Z.W. Birnbaum and W. Orlicz in 1930's \cite{bo}. This is a larger class of spaces which includes Lebesgue spaces, and has been used effectively in studying differential and integral equations. In harmonic analysis  the boundedness of operators on Orlicz spaces has been discussed widely. In \cite{hm}, Heinig and Maligranda provided four weight characterization for the Hardy operator on Orlicz spaces.  S. Bloom and R. Kerman \cite{bk} did characterize a four weight inequality on Orlicz spaces for generalized Hardy operator. Another characterization for the same was given by J. Q. Sun \cite{jqs2}.

Sawyer's duality principle, an important tool to study the boundedness of  linear operators defined on the class of monotone functions, initially given by E. Sawyer \cite{sa} for weighted Lebesgue spaces, was later developed for Orlicz spaces by Heinig and Kufner \cite{hk}. Also, they established a weight characterization for Hardy operator on the class of non-increasing functions in the framework of the  Orlicz spaces. In 1997, Drabek, Heinig and Kufner \cite{dhk} gave weight characterization of weighted modular inequalities for operators defined on monotone functions using the Sawyer's duality principle on Orlicz spaces given in \cite{hk}.   Our study extends the results of \cite{dhk}  and \cite{hk} for Dunkl-Hausdorff operator.
   
   In this section we provide some  definitions and preliminary results. In section 2, we give weight characterizations for Dunkl-Hausdorff operator on non-negative non-increasing function in terms of modular inequality. In Section 3, we prove some sufficient conditions for the boundedness of Dunkl-Hausdorff operator for non-negative non-increasing functions in terms of norm inequality. Lastly, we give a result establishing the equivalence of modular and norm inequalities under restricted conditions. 
 
\vspace{3pt}

We start with some definitions to be used throughout the paper.\\
Let $\varphi:[0,\infty) \rightarrow [0,\infty)$ be a function such that $\varphi(0)=0;~ \varphi(s)>0$ for $s>0;$ $\displaystyle \lim_{s\rightarrow\infty}\varphi(s) = \infty;~\varphi$ is right-continuous and non-decreasing on $[0,\infty).$ Then the function $\Phi$ defined as
\be \label{eqn1.1}
		\Phi(x):=\int_0^x \varphi(t)dt,~x \in [0,\infty)
\ee
is called an \emph{$N$-function} generated by $\varphi$.
 
\vspace{3pt}

\noindent An $N$-function $\Phi$ is said to satisfy the \emph{$\Delta_2$condition} (i.e., $\Phi \in \Delta_2$) if there exists a constant $\beta>0$ such that $\Phi(2t) \le \beta~ \Phi(t)$ for all $t > 0.$

\vspace{3pt}

\noindent By a weight, we shall mean a non-negative measurable function defined on an open set $\Omega \subseteq \mathbb R.$   Let $v$ be a weight function and  $\Phi \in \Delta_2,$ the \emph{weighted Orlicz space} denoted as $L_{\Phi(v)},$ is defined to be the collection of $v(x)dx-$measurable functions $f:\Omega \rightarrow \mathbb R$ such that
\[\rho_v(f,\Phi):=\int_\Omega \Phi(|f(x)|) v(x)dx \] is finite, where $\rho_v(f,\Phi)$ is termed as the \emph{modular} for $\Phi.$ The weighted Orlicz space is a Banach space with the Luxemburg norm, which is defined as	
	\[ \|f\|_{\Phi(v)}:= \inf \left\{\tau>0: \int_\Omega \Phi \lr \frac{|f(x)|}{\tau} \rr v(x) dx \le 1 \right\}.\]
	
\noindent If $\Phi$ is an $N$-function generated by  $\varphi,$ then the function defined as
\[\tilde{\Phi}(x):=\int_0^x \tilde{\varphi}(t)dt,~x \in [0,\infty)\]
is called the \emph{complementary function} to $\Phi,$ where $\tilde{\varphi}(t):=\displaystyle \sup_{\varphi(s) \le t}s.$ Sometimes, we refer to $\Phi$ and $\tilde{\Phi}$ as a pair of complementary functions.
Note that if the function $\varphi$  is strictly increasing then $\tilde{\varphi}$ is equal to $\varphi^{-1}.$

\vspace{5pt}

\noindent   For $\Omega = (0,\infty), $ following is the \emph{duality principle} in the weighted Orlicz space $L_{\Phi(v)}$:
    \be \label{b22}
     \sup\limits_{0\leq f} \frac{\int_{0}^{\infty}f(x)g(x)dx}{{\left\|f\right\|}_{\Phi(v)}} = \left\|\frac{g}{v}\right\|_{\tilde{\Phi}(v)}
    \ee
where $\tilde{\Phi}$ is the complementary function of the $N$-function $\Phi.$ For more on Orlicz spaces, see \cite{dhk, hk, kk,kr,pkjf,rr}.

Now onwards, throughout the paper, we assume  $\Omega = (0,\infty).$ All functions considered are non-negative and measurable, and all constants are positive. When $f$ is taken to be a non-negative non-increasing measurable function, we shall write it as $0 \le  {f\downarrow}.$ Similar notation will be used for non-decreasing case.  The notation  $\tilde{\Phi}$ will always denote the complementary function of $\Phi.$ The characteristic function of an interval with $0\leq a<b \leq \infty$ is denoted by $\chi_{(a,b)}.$ Additionally, we write  $A \approx B$ to indicate that the ratio $\frac{A}{B}$	is bounded  below and above by two constants $C_1$ and $C_2$ respectively.

\vspace{5pt} 

Analogous to the Sawyer's duality principle on weighted Lebesgue spaces \cite{sa}, Heinig and Kufner \cite{hk} established a duality principle for the cone of monotone functions in the framework of weighted Orlicz spaces. Below, we state this principle.

\vspace{3pt}
	
 \noindent {\bf Theorem A}\text{\bf \cite{hk}.}
  	Suppose $\Phi$ and $\tilde{\Phi}$ is a pair of complementary functions satisfying the $\Delta_{2}$condition.
 Let $Ig(x):=\int_{0}^{x} g, ~g \geq 0; ~Iv(x) := \int_{0}^{x} v$ with $(Iv)(\infty) = \infty,$ then 
    \be \label{b23}
  	\sup\limits_{0 \leq f \downarrow} \frac{\int_{0}^{\infty}f(x)g(x)dx}{{\left\|f\right\|}_{\Phi(v)}} \approx \left\|\frac{Ig}{Iv}\right\|_{\tilde{\Phi}(v)}.
  	\ee

   \begin{remark}  \label{rm1}
   For a linear operator \[(Tf)(x) := \int_{0}^{\infty} k(x,t) f(t) dt,~~ x>0,\] where $k$ is a non-negative kernel and $0\leq f\downarrow,$ the inequality 
	\[\left\| Tf \right\|_{\Psi(w)} \leq C \left\|f\right\|_{\Phi(v)},~0\leq f\downarrow\] 
	is equivalent to 
   \be \label{b25}
    \left\|\frac{I(T^{*}g)}{Iv}\right\|_{\tilde{\Phi}(v)} \leq C \left\| \frac{g}{w} \right\|_{\tilde{\Psi}(w)},   \hspace{2mm} g\geq 0,
   \ee
   where $T^{*}$ is the adjoint of $T,$  $\Psi$ and $\tilde{\Psi}$ are $N$-functions  and $v$ and $w$ are weight functions.  This equivalence, without proof is given in \cite{hk}, which can be proved by  using the arguments similar to that as given in \cite{sa}. For the sake of completion, below we give an explanation for the same.
	
	\vspace{5pt}
	
\noindent {\bf Explanation.} 
	Suppose the inequality (\ref{b25}) holds for any $g \ge 0.$ Let us consider $g'=hw$ for any $h\geq 0$ and $w$ be a weight function.  Clearly, $g' \ge 0.$ Using (\ref{b25}) with $g'$ we have 
	\be \label{b25'}
	\left\|\frac{I(T^{*}g')}{Iv}\right\|_{\tilde{\Phi}(v)} \leq C \left\| \frac{g'}{w} \right\|_{\tilde{\Psi}(w)}. 
	\ee	
	Then on applying Sawyer's duality principle \eqref{b23} for $T^*g'$ and using \eqref{b25'}, we obtain
  \begin{align}
  	\sup\limits_{0 \leq f \downarrow} \frac{\int_{0}^{\infty}f(x)(T^*g')(x)dx}{{\left\|f\right\|}_{\Phi(v)}} &\approx \left\|\frac{I(T^*g')}{Iv}\right\|_{\tilde{\Phi}(v)}\nonumber \\
  	  \sup\limits_{0 \leq f \downarrow} \frac{\int_{0}^{\infty}f(x)(T^*g')(x)dx}{{\left\|f\right\|}_{\Phi(v)}} & \leq C\left\| \frac{g'}{w} \right\|_{\tilde{\Psi}(w)} \nonumber \\
  	 \frac{\int_{0}^{\infty}f(x)(T^*g')(x)dx}{{\left\|f\right\|}_{\Phi(v)}} & \leq C\left\| \frac{g'}{w} \right\|_{\tilde{\Psi}(w)}, \hspace{2mm}  0\leq f\downarrow \nonumber \\  	  
     \frac{\int_{0}^{\infty}(Tf)(x)g'(x)dx}{{\left\|f\right\|}_{\Phi(v)}} & \leq C\left\| \frac{g'}{w} \right\|_{\tilde{\Psi}(w)}  \nonumber \\
  	 \frac{\int_{0}^{\infty}(Tf)(x)g'(x)dx}{\left\| \frac{g'}{w} \right\|_{\tilde{\Psi}(w)}} & \leq C {\left\|f\right\|}_{\Phi(v)}  \nonumber \\
  	  \frac{\int_{0}^{\infty}(Tf)(x)w(x)h(x) dx} {\left\|h\right\|_{\tilde{\Psi}(w)}} & \leq C {\left\|f\right\|}_{\Phi(v)}  \nonumber \\
  	  \sup\limits_{h \geq 0} \frac{\int_{0}^{\infty}(wTf)(x)h(x)dx}{\left\|h\right\|_{\tilde{\Psi}(w)}} & \leq C {\left\|f\right\|}_{\Phi(v)}  \nonumber \\	
  	  \left\| \frac{wTf}{w} \right\|_{\Psi(w)} & \leq C {\left\|f\right\|}_{\Phi(v)} \nonumber \\
  	  \left\| Tf \right\|_{\Psi(w)} & \leq C {\left\|f\right\|}_{\Phi(v)}, \hspace{2mm} 0\leq f\downarrow. \nonumber  
  \end{align}
    Likewise, the reverse part can be shown. 
	
 \end{remark}

 \noindent The modular version of the equivalence of norm inequalities mentioned in the Remark \ref{rm1} is proved in  \cite{dhk}. Precisely the following is known.

\vspace{3pt}

\noindent {\bf Theorem B}\text{\bf \cite{dhk}.} \emph{Suppose $\Phi, \tilde{\Phi} \in \Delta_{2}$  and $\Psi, \tilde{\Psi} $ are  $N$-functions, $w_0, ~ w_1, ~v$ are weight functions and $T$ is a positive linear operator defined on a collection of  monotone functions.}
   \emph{If $0 \leq f \downarrow$,$Ig(x)=\int_{0}^{x} g ~; Iv(x)=\int_{0}^{x} v $  and  $(Iv)(\infty) = \infty$,  then the modular inequalities}       
	\be \label{b27}
   	   {\Psi}^{-1}\lge \int_{0}^{\infty} \Psi \lek w_{1}(x) (Tf)(x) \rek w_{0}(x) dx \rge  \leq {\ \Phi}^{-1}\lge \int_{0}^{\infty} \Phi \lek  Cf(x) \rek v(x) dx \rge 
      \ee     
      \emph{and}      
      \[{\Psi}^{-1}\lge \int_{0}^{\infty} \Psi \lek w_{1}(x) T(I^{*}g)(x) \rek w_{0}(x) dx \rge  \leq {\Phi}^{-1}\lge \int_{0}^{\infty} \Phi \lek  C\frac{(Iv)(x) g(x)}{v(x)} \rek v(x) dx \rge \]
 \emph{with functions $g \geq 0$ are equivalent, where $I^*g(x):=\int_{x}^{\infty} g.$}
     
 \vspace{3pt}    
		
Now, before mentioning some results related to the modular integral inequalities for operators of Hardy type, we define the order relation in between two $N$-functions.
 
\noindent For $N$-functions $\Phi$ and $\Psi,$ we say that \emph{$\Phi$ is less than $\Psi,$} denoted by $\Phi \prec \Psi,$ if there is a constant $C>0,$ such that   		
  		\[\sum_{i} \Psi \circ \Phi^{-1}(a_{i}) \leq C ~\Psi \circ \Phi^{-1}\lr \sum_{i} a_{i} \rr \]  		
for every non-negative sequence $\{{a_{i}}\}.$ 
If $\Psi \circ \Phi^{-1}$ is convex, then $\Phi \prec \Psi.$

In \cite{bk},  authors have proved the modular integral inequalities for operators of Hardy type by assuming the condition $\Psi \circ \Phi^{-1}$ is convex but in \cite{jqs2} these inequalities have been proved with the weaker condition $\Phi \prec \Psi.$ Below we state the results given in \cite{jqs2}. 

\noindent An operator of the form 
 \[(Kf)(x) := \int_{0}^{x} k(x,t) f(t) dt \]
 is known as generalized Hardy operator (GHO), where
      \begin{enumerate}[(a)]
      	\item $k: \{(x,t)\in \mathbb{R}^{2}: 0<t<x<\infty\} \rightarrow [0,\infty);$
      	\item $k(x,t)\geq 0$ is non-decreasing ($\uparrow$) in $x,$  non-increasing ($\downarrow$) in $t;$ and
      	\item $k(x,y) \leq D(k(x,t)+ k(t,y))$ whenever $0\leq y\leq t\leq x <\infty$ for some constant D.
      \end{enumerate}

\noindent {\bf Theorem C}\text{\bf\cite{jqs2}.}  
	\emph {Let $\Phi$ and $\Psi$ be $N$-functions such that $\Phi \prec \Psi,$ and $K$ be  GHO. Let $a,b,v$ and $w$ be weight functions. Then there exists a constant $A$ such that}
	
	\be \label{b21}
{\Psi}^{-1}\lge \int_{0}^{\infty} \Psi \lek a(x) (Kf)(x) \rek w(x) dx \rge \leq {\Phi}^{-1}\lge \int_{0}^{\infty} \Phi \lek  A b(x) f(x) \rek v(x) dx \rge 
\ee
\emph {holds for all functions f if and only if there exists a constant $C$ such that both of the following inequalities}
	\[   
	{\Psi}^{-1}\lge \int_{r}^{\infty} \Psi \lek \frac{a(x)}{C} \left\| \frac{k(r,.)\chi_{(0,r)}} {\varepsilon vb} \right\|_{\tilde{\Phi}(\varepsilon v)} \rek w(x) dx \rge \leq {\Phi}^{-1}(1/\varepsilon)   
	\]
	\emph {and}
	\[
	{\Psi}^{-1}\lge \int_{r}^{\infty} \Psi \lek \frac{a(x)}{C} \left\| \frac{\chi_{(0,r)}} {\varepsilon vb} \right\|_{\tilde{\Phi}(\varepsilon v)} k(x,r) \rek w(x) dx \rge \leq {\Phi}^{-1}(1/\varepsilon)
	\]
	\emph {hold for all $\varepsilon >0, r>0.$}
	
\vspace{3pt}

\noindent For the adjoint $K^{*}$ of the  above operator $K,$ defined as 
\[   (K^{*}g)(x) := \int_{x}^{\infty} k(t,x) g(t) dt,\]
we have the following result:

\vspace{3pt}

\noindent {\bf Theorem D}\text{\bf \cite{jqs2}.} 
	\emph{Let $\Phi$ and $\Psi$ be $N$-functions such that $\Phi \prec \Psi,$ and $K^{*}$ be the adjoint of $GHO.$ Let $a,b,v$ and $w$ be weight functions. Then there exists a constant $A$ such that the inequality}
	\[
	{\Psi}^{-1}\lge \int_0^\infty \Psi \lek a(x) (K^{*}g)(x) \rek w(x) dx \rge \leq {\Phi}^{-1}\lge \int_{0}^{\infty} \Phi \lek  A b(x) g(x) \rek v(x) dx \rge 
	\]   	
	\emph{holds for all functions $g$ if and only if there exists a constant $C$ such that both of the following inequalities} 
	\[   
	\Psi^{-1}\lge \int_{0}^{r} \Psi \lek \frac{a(x)}{C} \left\| \frac{k(.,r)\chi_{(r,\infty)}} {\varepsilon vb} \right\|_{\tilde{\Phi}(\varepsilon v)} \rek w(x) dx \rge \leq {\Phi}^{-1}(1/\varepsilon)   
	\]
	\emph{and}
	\[
	{\Psi}^{-1}\lge \int_{0}^{r} \Psi \lek \frac{a(x)}{C} \left\| \frac{\chi_{(r,\infty)}} {\varepsilon vb} \right\|_{\tilde{\Phi}(\varepsilon v)} k(r,x) \rek w(x) dx \rge \leq {\Phi}^{-1}(1/\varepsilon)
	\]
	\emph{hold for all} $\varepsilon >0, r>0.$
	
	\vspace{3pt}
   
 We shall require the following result from \cite{dhk}, which gives the characterization of the  weighted modular inequality for a linear operator $T$ defined on $\sigma$-finite measure spaces, in terms of a certain weighted  Orlicz-Luxemburg norm inequality. 

\vspace{3pt}
 
\noindent {\bf Proposition E}\text{\bf \cite{dhk}.} 
  	\emph{Suppose $(X,d\mu), (Y,d\nu)$ are $\sigma$-finite measure spaces and $T$ is a linear operator mapping measurable functions on $X$ to measurable functions on $Y.$ If $\Phi$ and $\Psi$ are $N$-functions, then the modular inequality
  	 \ben \label{b28}
  	{\Psi}^{-1}\lge \int_{Y} \Psi \lek w(y) |(Tf)(y)| \rek d\nu(y) \rge \leq {\
  		\Phi}^{-1}\lge \int_{X} \Phi \lek  Cu(x)|f(x)| \rek d\mu(x) \rge 
  	\een
  	is satisfied if and only if for every $\varepsilon>0,$ the inequality
  	\ben \label{b29}
  	\left\| wTf \right\|_{\Psi(\varepsilon_{\Psi}d\nu)} \leq C \left\| uf \right\|_{\Phi(\varepsilon_{\Phi}d\mu)} 
  	\een
  	holds, where $\varepsilon_{\Psi} = 1/\Psi(1/\varepsilon), \varepsilon_{\Phi} = 1/\Phi(1/\varepsilon).$}
 
\vspace{3pt}
		
		Finally, we mention that our work is motivated by Sawyer’s duality principle \cite{sa} for weighted Lebesgue spaces, which was later extended to weighted Orlicz spaces by Heinig and Kufner in \cite{hk}. As an application of this principle, we have proved a weight characterization required for the modular and norm inequalities to hold for the Dunkl-Hausdorff operator on non-negative, non-increasing functions in the framework of weighted Orlicz spaces.
		
   \section{Modular Inequalities}
   
   We start with the following definition.\\
	
	Let $\phi \in L^{1}_{loc}(\mathbb{R}),$ then the Dunkl-Hausdorff operator $\mathcal{H}_{\alpha,\phi}$ for $\alpha \in \mathbb{R}$ is defined as 
   \ben
   (\mathcal{H}_{\alpha,\phi} f)(x) := \int_{0}^{\infty} \frac{\phi(t)}{t^{2\alpha + 2}} f\lr\frac{x}{t}\rr dt.
   \een
   
   \noindent For $\alpha = -1/2,$ the operator $\mathcal{H}_{\alpha,\phi}$ becomes the Hausdorff operator 
   \ben
   (\mathcal{H}_{\phi} f)(x) := \int_{0}^{\infty} \frac{\phi(t)}{t} f\lr\frac{x}{t}\rr dt.
   \een   
   The Hausdorff operator returns different operators for suitable values of $\phi.$ For example, \\
the Hardy averaging operator
\[(\mathcal{A}f)(x) := \frac{1}{x}\int_{0}^{x} f(t) dt\] 
for $\phi(t) = \frac{\chi_{(1,\infty)}(t)}{t};$ the adjoint of Hardy averaging operator 
 \[(\mathcal{A}^{*}f)(x) := \int_{x}^{\infty} \frac{f(t)}{t} dt\] \\
for $\phi(t) = \chi_{(0,1)}(t).$ For more on Dunkl-Hausdorff and Hausdorff operator, see \cite{bg,bs,bs2,cfw,jfj}.\\

Below we give our first main result.  
	
   \begin{theorem} \label{t1}
   	Suppose $\Phi$ and $\Psi$ be N-functions such that $\Phi, \tilde{\Phi} \in \Delta_{2}$ and $\Phi \prec \Psi.$ Let $(Iv)(\infty) = \infty$ and for $\alpha \in \mathbb{R},$ we have
   		\begin{eqnarray}\label{a21}
   	      \begin{cases} 
   			C_{1} \leq \frac{\phi(x)}{x^{2\alpha}} \leq C_{2} , &  x \geq 1 \\
   			C_{3} \leq  \frac{\phi(x)}{x^{(2\alpha +1)}} \leq C_{4} , & 0< x < 1
   		\end{cases}
   	\end{eqnarray}
   for some constants $C_{i}, ~i = 1,2,3,4.$ Then the modular inequality
   \be\label{a22}
   	{\Psi}^{-1}\lge \int_{0}^{\infty} \Psi \lek w_{1}(x) (\mathcal{H}_{\alpha,\phi} f)(x) \rek w_{0}(x) dx \rge \leq {\Phi}^{-1}\lge \int_{0}^{\infty} \Phi \lek  Cf(x) \rek v(x) dx \rge 
   \ee   
  holds for all $0\leq f \downarrow$ if and only if there are constants $A_1,A_2,A_3$ such that for all $r>0, ~\varepsilon>0,$ the following inequalities hold:\\
    \begin{align}\label{a23}
  	{\Psi}^{-1}\lge \int_{r}^{\infty} \Psi \lek \frac{w_{1}(x)}{A_1x} \left\| \frac{x\chi_{(0,r)}}{\varepsilon(Iv)} \right\|_{\tilde{\Phi}(\varepsilon v)} \rek w_{0}(x) dx \rge \leq {\Phi}^{-1}(1/\varepsilon)
  	\end{align}
  	\begin{align} \label{a24}  		
  	{\Psi}^{-1}\lge \int_{0}^{r} \Psi \lek \frac{w_{1}(x)}{A_2} \left\| \frac{\chi_{(r,\infty)}}{\varepsilon(Iv)} \right\|_{\tilde{\Phi}(\varepsilon v)} \rek w_{0}(x) dx \rge \leq {\Phi}^{-1}(1/\varepsilon)
  	\end{align}
  	\begin{align}\label{a25}
  	{\Psi}^{-1}\lge \int_{0}^{r} \Psi \lek \frac{w_{1}(x)}{A_3} \left\| \frac{\chi_{(r,\infty)}\ln\lr\frac{.}{r}\rr} {\varepsilon(Iv)} \right\|_{\tilde{\Phi}(\varepsilon v)} \rek w_{0}(x) dx \rge \leq {\Phi}^{-1}(1/\varepsilon) 
    \end{align}
   	\begin{align} \label{a26}
    {\Psi}^{-1}\lge \int_{0}^{r} \Psi \lek \frac{w_{1}(x)}{A_3} \left\| \frac{\chi_{(r,\infty)}} {\varepsilon(Iv)} \right\|_{\tilde{\Phi}(\varepsilon v)} \ln\lr\frac{r}{x}\rr \rek w_{0}(x) dx \rge \leq {\Phi}^{-1}(1/\varepsilon)
   \end{align}
   \end{theorem}
  
   \proof
  On using Theorem B, we have that the inequality \eqref{a22} holds for $0 \le f \downarrow$ if and only if the following inequality 
   \be \label{a27}
   	{\Psi}^{-1}\lge \int_{0}^{\infty}  \Psi \lek w_{1}(x)  \mathcal{H}_{\alpha,\phi}(I^{*}h)(x) \rek w_{0}(x) dx \rge 
		 \leq {\Phi}^{-1}\lge \int_{0}^{\infty} \Phi \lek  C\frac{(Iv)(x) h(x)}{v(x)} \rek v(x) dx \rge 
   \ee 
  holds for the functions $h\geq0.$
   
   \noindent On computing $\mathcal{H}_{\alpha,\phi}(I^{*}h),$ we have 
   \begin{align} \label{a28}
   \mathcal{H}_{\alpha,\phi}(I^{*}h)(x) & = \int_{0}^{\infty} \frac{\phi(x/y)}{x^{2\alpha+1}} y^{2\alpha} (I^{*}h)(y) dy \nonumber \\
   &= \int_{0}^{\infty} \frac{\phi(x/y)}{x^{2\alpha+1}} y^{2\alpha} \lr\int_{y}^{\infty} h(t) dt \rr dy \nonumber \\
   &= \int_{0}^{x} \frac{\phi(x/y)}{x^{2\alpha+1}} y^{2\alpha} \lr\int_{y}^{\infty} h(t) dt \rr dy + \int_{x}^{\infty} \frac{\phi(x/y)}{x^{2\alpha+1}} y^{2\alpha} \lr\int_{y}^{\infty} h(t) dt \rr dy \nonumber \\
   &=: \mathcal{I}_{1} + \mathcal{I}_{2}.
   \end{align}
   
  \noindent On using \eqref{a21} and performing change of variable, we estimate $\mathcal{I}_{1}$ as below
  \begin{align}\label{a29}
   \mathcal{I}_{1} &= \int_{0}^{x} \frac{\phi(x/y)}{x^{2\alpha+1}} y^{2\alpha} \lr\int_{y}^{\infty} h(t) dt \rr dy \nonumber \\
   &= \int_{0}^{x} \frac{\phi(x/y)}{x^{2\alpha+1}} y^{2\alpha} \lr\int_{y}^{x} h(t) dt \rr dy + \int_{0}^{x} \frac{\phi(x/y)}{x^{2\alpha+1}} y^{2\alpha} \lr\int_{x}^{\infty} h(t) dt \rr dy  \nonumber \\
   &\leq C_2 \lek \int_{0}^{x} \frac{x^{2\alpha}}{y^{2\alpha}} \frac{y^{2\alpha}}{x^{2\alpha+1}}  \lr\int_{y}^{x} h(t) dt \rr dy + \int_{0}^{x} \frac{x^{2\alpha}}{y^{2\alpha}} \frac{y^{2\alpha}}{x^{2\alpha+1}} \lr\int_{x}^{\infty} h(t) dt \rr dy \rek \nonumber \\
   &= C_2 \lek \int_{0}^{x} \frac{1}{x} \lr\int_{y}^{x} h(t) dt \rr dy + \int_{0}^{x} \frac{1}{x} \lr\int_{x}^{\infty} h(t) dt \rr dy \rek \nonumber \\
   & = C_2 \lek \frac{1}{x} \int_{0}^{x}  t h(t) dt  +  \int_{x}^{\infty} h(t) dt \rek. 
  \end{align}
    
   \noindent Similarly, using \eqref{a21} we obtain
   \begin{eqnarray} \label{a210}
   \mathcal{I}_1\geq C_{1} \lek \frac{1}{x} \int_{0}^{x}  t h(t) dt  +  \int_{x}^{\infty} h(t) dt \rek.
   \end{eqnarray} 
   On combining \eqref{a29} and \eqref{a210} we have
    \begin{eqnarray} \label{a211}
     \mathcal{I}_1 \approx \frac{1}{x} \int_{0}^{x}  t h(t) dt  +  \int_{x}^{\infty} h(t) dt.  
    \end{eqnarray}
   
    \noindent Likewise, using \eqref{a21} we estimate $\mathcal{I}_{2}$ as below
    \begin{align} \label{a212}
     \mathcal{I}_{2} &= \int_{x}^{\infty} \frac{\phi(x/y)}{x^{2\alpha+1}} y^{2\alpha} \lr\int_{y}^{\infty} h(t) dt \rr dy \nonumber \\
     &\leq C_4 \int_{x}^{\infty} \frac{x^{2\alpha+1}}{y^{2\alpha+1}} \frac{y^{2\alpha}}{x^{2\alpha+1}}  \lr\int_{y}^{\infty} h(t) dt \rr dy \nonumber \\
     &= C_{4} \int_{x}^{\infty} \frac{1}{y} \lr\int_{y}^{\infty} h(t) dt \rr dy \nonumber \\
     &= C_{4} \int_{x}^{\infty} \ln\lr\frac{t}{x}\rr  h(t) dt.
    \end{align}
   
  \noindent Also,  using \eqref{a21}, we have
   \begin{eqnarray}\label{a213}
   	\mathcal{I}_{2}\geq C_{3} \int_{x}^{\infty} \ln\lr\frac{t}{x}\rr  h(t) dt. 
   \end{eqnarray}

   \noindent So, \eqref{a212} and \eqref{a213} combined together give
    \begin{eqnarray}\label{a214}
    \mathcal{I}_{2} \approx \int_{x}^{\infty} \ln\lr\frac{t}{x}\rr  h(t) dt. 
    \end{eqnarray}

   \noindent Consequently, \eqref{a27} is equivalent to    
   	\begin{align*}   	   	
   	{\Psi}^{-1}\lge \int_{0}^{\infty} \Psi \lek w_{1}(x) \lr \frac{1}{x} \int_{0}^{x}  t h(t) dt   +  \int_{x}^{\infty} h(t) dt + \int_{x}^{\infty} \ln\lr\frac{t}{x}\rr  h(t) dt  \rr  \rek w_{0}(x) dx \rge \\
      \leq {\Phi}^{-1}\lge \int_{0}^{\infty} \Phi \lek  C\frac{(Iv)(x) h(x)}{v(x)} \rek v(x) dx \rge.      		
   	\end{align*}   
   
   \noindent Since $\Phi \in \Delta_{2},$ using the convexity of $\Psi$ and $\Phi,$ it is easily observed that the above inequality is equivalent to the following three inequalities:
   
   \begin{multline}\label{a215}
   	{\Psi}^{-1}\lge \int_{0}^{\infty} \Psi \lek w_{1}(x) \lr \frac{1}{x} \int_{0}^{x}  t h(t) dt  \rr  \rek w_{0}(x) dx \rge
   	\\ \leq {\Phi}^{-1}\lge \int_{0}^{\infty} \Phi \lek  C\frac{(Iv)(x) h(x)}{v(x)} \rek v(x) dx \rge, 	
   \end{multline}

    \begin{multline}\label{a216}
    {\Psi}^{-1}\lge \int_{0}^{\infty} \Psi \lek w_{1}(x) \lr   \int_{x}^{\infty} h(t) dt \rr \rek w_{0}(x) dx \rge \\ 
    \leq {\Phi}^{-1}\lge \int_{0}^{\infty} \Phi \lek  C\frac{(Iv)(x) h(x)}{v(x)} \rek v(x) dx \rge
    \end{multline}
and   
   \begin{multline}\label{a217}
   {\Psi}^{-1}\lge \int_{0}^{\infty} \Psi \lek w_{1}(x)  \lr  \int_{x}^{\infty} \ln\lr\frac{t}{x}\rr  h(t) dt  \rr  \rek w_{0}(x) dx \rge \\ \leq {\Phi}^{-1}\lge \int_{0}^{\infty} \Phi \lek  C\frac{(Iv)(x) h(x)}{v(x)} \rek v(x) dx \rge. 
   \end{multline}

   \noindent Now, the operator in \eqref{a215} is the Hardy averaging operator, and thus, on using Theorem C for $k(x,t)=1,$ \eqref{a215} holds if and only if \eqref{a23} holds. Likewise on using Theorem D with $k(x,t)=1,$ \eqref{a216} holds if and only if \eqref{a24} holds. Lastly, for $k(t,x)=\ln \lr\frac{t}{x}\rr$ being taken in Theorem D, the inequality \eqref{a217} holds if and only if \eqref{a25} and \eqref{a26} hold.  $ \hfill \square$ \\
   
	\begin{remark}
	In Theorem 2.1, under the condition (2.1) on $\phi,$ the Dunkl-Hausdroff operator is equivalent to the  Calderan operator (Hardy-Littlewood-Polya operator), which is defined as 
\be \label{eqn218}
Ch(x):= \frac{1}{x}\int_0^x h(t)dt+ \int_x^\infty\frac{h(t)}{t} dt,
\ee
for non-negative measurable functions $ h \ge 0.$ I.e., 
${H}_{\alpha, \phi}h \approx {C}h$ for $h \ge 0.$  So, under the condition (2.1) on $\phi,$ the  study of the boundedness of the Dunkl-Hausdroff operator for the class of non-negative functions is equivalent to the study of  the boundedness of the Hardy-Littlewood-Polya operator. But, to study the boundedness of the Dunkl-Hausdroff operator on the class of non-negative  non-increasing functions, after applying Theorem B [9], the computations are done on ${H}_{\alpha, \phi}(I^* h),$ where $I^* h(y)= \int_y^\infty h(t) dt.$ Under  the condition (2.1) on $\phi,$ it is shown in the proof of Theorem 2.1 that ${H}_{\alpha, \phi}(I^* h)$ is equivalent to the following:
\be
\frac{1}{x}\int_0^x  t h(t)dt+ \int_x^\infty h(t) dt+\int_x^\infty \ln \lr\frac{t}{x}\rr h(t) dt, \nonumber
\ee
by (2.8), (2.11) and (2.14). I.e., for $ h \ge 0,$ we have 
\[{H}_{\alpha, \phi}(I^* h) \approx \frac{1}{x}\int_0^x  t h(t)dt+ \int_x^\infty h(t) dt+\int_x^\infty \ln \lr\frac{t}{x}\rr h(t) dt\]
which is clearly different from ${C}h$ in (\ref{eqn218}).\\
\end{remark}

	A result similar to Theorem 2.1 can also be proved for quasi Dunkl-Hausdorff operator $\mathcal{H}_{\alpha,\phi}^{*},$ which is the adjoint of the operator $\mathcal{H}_{\alpha,\phi}$ defined as 
  \[
  \left(\mathcal{H}_{\alpha,\phi}^{*} g\right)(x) := x^{2\alpha} \int_{0}^{\infty} \frac{\phi(y/x)}{y^{(2\alpha+1)}} g(y) dy. 
  \]
	To that we have the following.
 \begin{theorem}
  Suppose $\Phi$ and $\Psi$ are N- functions such that $\Phi, \tilde{\Phi} \in \Delta_{2}$ and $\Phi \prec \Psi.$ Let $(Iv)(\infty) = \infty$ and for $\alpha \in \mathbb{R},$ the condition \eqref{a21} holds for some constants $C_{i}, ~i = 1,2,3,4.$ Then the inequality
  \be\label{a218}
  	{\Psi}^{-1}\lge \int_{0}^{\infty} \Psi \lek w_{1}(x) (\mathcal{H}_{\alpha,\phi}^{*} f)(x) \rek w_{0}(x) dx \rge \leq {\Phi}^{-1}\lge \int_{0}^{\infty} \Phi \lek  Cf(x) \rek v(x) dx \rge 
  \ee
	holds for $0\leq f \downarrow$ if and only if \eqref{a23}, \eqref{a24}, \eqref{a25} and \eqref{a26} hold.
  \end{theorem}
   
   \proof
  Similar to the proof of the above theorem, on using Theorem B, the inequality \eqref{a218} turns out to be equivalent to the following inequality:
 \ben
  {\Psi}^{-1}\lge \int_{0}^{\infty} \Psi \lek w_{1}(x) \mathcal{H}_{\alpha,\phi}^{*}(I^{*}h)(x) \rek w_{0}(x) dx \rge  \leq {\Phi}^{-1}\lge \int_{0}^{\infty} \Phi \lek  C\frac{(Iv)(x) h(x)}{v(x)} \rek v(x) dx \rge 
  \een
  for the functions $h\geq0.$  
   
  \noindent On computing $\mathcal{H}_{\alpha,\phi}^{*}(I^{*}h),$ we get
  \begin{align*} 
  	\mathcal{H}_{\alpha,\phi}^{*}(I^{*}h)(x) &= \int_{0}^{\infty} \frac{\phi(y/x)}{y^{2\alpha+1}} x^{2\alpha} I^{*}h(y) dy \nonumber \\
  	&= \int_{0}^{\infty} \frac{\phi(y/x)}{y^{2\alpha+1}} x^{2\alpha} \lr\int_{y}^{\infty} h(t) dt \rr dy \nonumber \\
  	&= \int_{0}^{x} \frac{\phi(y/x)}{y^{2\alpha+1}} x^{2\alpha} \lr\int_{y}^{\infty} h(t) dt \rr dy + \int_{x}^{\infty} \frac{\phi(y/x)}{y^{2\alpha+1}} x^{2\alpha} \lr\int_{y}^{\infty} h(t) dt \rr dy \nonumber \\
  	&=: \mathcal{I}_{3} + \mathcal{I}_{4}.
  \end{align*}
   
   \noindent Then on using \eqref{a21} it is easy to verify that
   \begin{align*}
   \mathcal{I}_3 \approx \frac{1}{x} \int_{0}^{x}  t h(t) dt  +  \int_{x}^{\infty} h(t) dt  
   \end{align*}
    and 
   \begin{align*}
   	\mathcal{I}_4 \approx \int_{x}^{\infty} \ln\lr\frac{t}{x}\rr  h(t) dt.  
   \end{align*}
   Now the proof follows on using Theorems C and D with suitable values of the kernel $k(\cdot,\cdot).$ \hfill $\square$
		
		\vspace{3pt}
   
 From the above, we may deduce the following corollaries for the Hausdorff operator $\mathcal{H}_{\phi}$ and its adjoint $\mathcal{H}^*_{\phi}$ respectively. 
\begin{corollary} \label{c1}
   	 Suppose $\Phi$ and $\Psi$ are N- functions such that $\Phi, \tilde{\Phi} \in \Delta_{2}$ and $\Phi \prec \Psi.$ Let $(Iv)(\infty) = \infty$ and  we have
   	\begin{eqnarray} \label{a219} 
   		\begin{cases} 
   			C_{1} \leq x{\phi(x)} \leq C_{2} , &  x \geq 1 \\
   			C_{3} \leq {\phi(x)} \leq C_{4} , &  0< x < 1
   		\end{cases}
   	\end{eqnarray}
   	for some constants $C_{i}, ~i = 1,2,3,4.$ Then the inequality
   	\begin{align*}
   		{\Psi}^{-1}\lge \int_{0}^{\infty} \Psi \lek w_{1}(x) (\mathcal{H}_{\phi}f)(x) \rek w_{0}(x) dx \rge \leq {\Phi}^{-1}\lge \int_{0}^{\infty} \Phi \lek  Cf(x) \rek v(x) dx \rge 
   	\end{align*}   
   	holds for $0\leq f \downarrow$ if and only if \eqref{a23}, \eqref{a24}, \eqref{a25} and \eqref{a26} hold.
   \end{corollary}

   \begin{corollary}
   	 Suppose $\Phi$ and $\Psi$ are N- functions such that $\Phi, \tilde{\Phi} \in \Delta_{2}$ and $\Phi \prec \Psi.$ Let $(Iv)(\infty) = \infty$ and \eqref{a219} holds for some constants $C_{i}, ~i = 1,2,3,4.$ Then the inequality
   	\begin{align*}
   		{\Psi}^{-1}\lge \int_{0}^{\infty} \Psi \lek w_{1}(x) (\mathcal{H}_{\phi}^{*}f)(x) \rek w_{0}(x) dx \rge \leq {\Phi}^{-1}\lge \int_{0}^{\infty} \Phi \lek  Cf(x) \rek v(x) dx \rge 
   	\end{align*}   
   	holds for $0\leq f \downarrow$ if and only if \eqref{a23}, \eqref{a24}, \eqref{a25} and \eqref{a26} hold.
   \end{corollary}

   \begin{remark}
   	On taking $\Phi(t) = \frac{t^{p}}{p}$ and $\Psi(t) = \frac{t^{q}}{q}, 1<p\leq q<\infty,$ the Theorem \ref{t1} reduces to a result for Dunkl-Hausdorff operator when considered on non-increasing functions in weighted Lebesgue spaces. An $n$-dimensional analogue of the later has been given by the authors in \cite{jjsm} for the Dunkl-Hausdorff operator on non-negative radially decreasing functions. Also, if we consider the Corollary \ref{c1} for the above taken values of $\Phi$  and $\Psi,$ we get the result for Hausdorff operator for  non-increasing functions on weighted Lebesgue spaces as given  in \cite{jj}.
   \end{remark}

\vspace{3pt}

   \noindent Consider the function
   \begin{eqnarray*}
   	\phi(t) = 
   	\begin{cases}
   		\frac{1}{t} , &  t \geq 1 \\   	
   		1 , &  0<t < 1.
   	\end{cases} \\
   \end{eqnarray*}
   
   \noindent Note that $\phi$ satisfies the conditions of Corollary \ref{c1}, and the operator $\mathcal{H}_{\phi}$ reduces to the operator $\mathcal{A} + \mathcal{A}^{*},$ a Calder\'{o}n type operator. Thus, we have the following:
   \begin{corollary}
   	 Suppose $\Phi$ and $\Psi$ are N- functions such that $\Phi, \tilde{\Phi} \in \Delta_{2}$ and $\Phi \prec \Psi.$ Let $(Iv)(\infty) = \infty.$ Then the inequality
   	 \begin{align*}
   	 	{\Psi}^{-1}\lge \int_{0}^{\infty} \Psi \lek w_{1}(x) (\mathcal{A}+\mathcal{A}^{*})(f)(x) \rek w_{0}(x) dx \rge \leq {\Phi}^{-1}\lge \int_{0}^{\infty} \Phi \lek  Cf(x) \rek v(x) dx \rge 
   	 \end{align*}   
   	 holds for $0\leq f \downarrow$ if and only if \eqref{a23}, \eqref{a24}, \eqref{a25} and \eqref{a26} hold.
   \end{corollary}  
       
		On taking a restriction on the support of the function $\phi,$  we have the following lemma:
  
  \begin{lemma}
   Let $\phi \in {L^{1}_{loc}}(\mathbb{R}^{+})$ be such that  $supp (\phi) \subseteq [1, \infty),$ and $C_{1} \leq t \phi(t) \leq C_{2}, ~ t \geq 1$ for constants $C_1, C_2 >0.$ Then $\mathcal{H}_\phi \approx \mathcal{A}.$ 
  \end{lemma}		
  \proof 
  Using the fact that supp$\phi\subseteq [1, \infty)$ and change of variables, we get that
  \begin{align*}
   (\mathcal{H}_\phi f)(x) &= \int_{0}^{\infty} \frac{\phi(y)}{y} f\left({\frac{x}{y}} \right) dy \\
   &=\int_{1}^{\infty} \frac{\phi(y)}{y} f\left({\frac{x}{y}} \right) dy \\
   &=\int_{0}^{x} \frac{\phi(x/t)}{t} f(t) dt.  
  \end{align*}
  
  \noindent Now using that $C_{1} \leq t \phi(t) \leq C_{2}, ~ t \geq 1$ for constants $C_1, C_2 >0,$ we obtain
  \[
  \int_{0}^{x} \frac{\phi(x/t)}{t} f(t) dt \approx \frac{1}{x} \int_{0}^{x} f(t) dt.
  \]
  Thus $\mathcal{H}_\phi f(x) \approx \mathcal{A} f(x),~ x\in [0,\infty).$   \hfill $\square$ 
  	
		\vspace{3pt}
		
		Consequently, we have the following.
		
    \begin{corollary} \label{c25}
   	Suppose $\Phi$ and $\Psi$ are N- functions such that $\Phi, \tilde{\Phi} \in \Delta_{2}$ and $\Phi \prec \Psi.$ Let $(Iv)(\infty) = \infty, ~supp(\phi) \subseteq [1,\infty)$ and 
   	$C_{1} \leq x{\phi(x)} \leq C_{2} , x \geq 1$ for some constant $C_{i}, ~i = 1,2.$ Then the following are equivalent: 
   	\begin{enumerate} [(i)]
   		\item 	The inequality
   		\be
   		{\Psi}^{-1}\lge \int_{0}^{\infty} \Psi \lek w_{1}(x) (\mathcal{H}_{\phi}f)(x) \rek w_{0}(x) dx \rge \leq {\Phi}^{-1}\lge \int_{0}^{\infty} \Phi \lek  Cf(x) \rek v(x) dx \rge
   		\ee
   		holds for $0\leq f \downarrow.$
   		\item The inequality
   		\be \label{c29}
   		{\Psi}^{-1}\lge \int_{0}^{\infty} \Psi \lek w_{1}(x) (\mathcal{A}f)(x) \rek w_{0}(x) dx \rge \leq {\Phi}^{-1}\lge \int_{0}^{\infty} \Phi \lek  Cf(x) \rek v(x) dx \rge
   		\ee
   		holds for $0\leq f \downarrow.$
   		\item The conditions \eqref{a23} and \eqref{a24} hold.
   	\end{enumerate}
   \end{corollary} 
   \proof The above equivalence follows by showing that $(ii)$ holds $\iff (i)$ holds $\iff (iii)$ holds.   \hfill $\square$
   
   \begin{remark}
   	In the above corollary, we have given the weight characterization for modular inequalities in Orlicz spaces for Hardy averaging operator $\mathcal{A}$ on non-negative non-increasing functions. A similar result for operator $\mathcal{A}$ has been given in \cite{jqs2} (Theorem 4.3, given below). Therefore, using Corollary 2.1 of \cite{dhk} given below, we may say that the weight characterization in Corollary \ref{c25} is equivalent to the weight characterization given in (\cite{jqs2}, Theorem 4.3). 	
   \end{remark}
	
	\vspace{3pt}

 \noindent{\bf Lemma F}{\bf (\cite{jqs2}, Theorem 4.3).}
  \emph{Suppose $\Phi$ and $\Psi$ are N- functions such that $\Phi, \tilde{\Phi} \in \Delta_{2}$ and $\Phi \prec \Psi.$ Then the inequality \eqref{c29} with $w_1=1,$ holds for all $0\leq f\downarrow$ if and only if there is a constant B such that
  \begin{align*}
  \Psi^{-1}\lge\Psi(\varepsilon)(Iw_0)(r)\rge \leq \Phi^{-1}\lge \Phi(B\varepsilon) (Iv)(r)\rge
  \end{align*}
and 
  \begin{align*}
  	{\Psi}^{-1}\lge \int_{r}^{\infty} \Psi \lek \frac{1}{Bx} \left\| \frac{x\chi_{(0,r)}}{\varepsilon(Iv)} \right\|_{\tilde{\Phi}(\varepsilon v)} \rek w_{0}(x) dx \rge \leq {\Phi}^{-1}(1/\varepsilon)
  \end{align*}
   hold for all $\varepsilon, r >0.$}

\vspace{3pt}

 \noindent{\bf Lemma G}{\bf (\cite{dhk}, Corollary 2.1).}
  \emph{Suppose $\Phi$ and $\Psi$ are N- functions such that $\Phi \prec \Psi.$ Let $\Phi, \tilde{\Phi} \in \Delta_{2}$ and $(Iv)(\infty) = \infty.$ If  $0\leq f\downarrow,$ then the following are equivalent:
  \begin{enumerate} [(i)]
  \item 	For weights $w_0, w_1, v$ and a constant $C>0$ the inequality
	\ben
  	{\Psi}^{-1}\lge \int_{0}^{\infty} \Psi \lek w_{1}(x) f(x) \rek w_{0}(x) dx \rge \leq {\Phi}^{-1}\lge \int_{0}^{\infty} \Phi \lek  Cf(x) \rek v(x) dx \rge
    \een
		holds.
    \item There exists a constant $B>0,$ such that for all $r>0, \varepsilon>0,$ the estimate
    \ben
   	{\Psi}^{-1}\lge \int_{0}^{r} \Psi \lek \frac{1}{B} \left\| \frac{\chi_{(r,\infty)}}{\varepsilon(Iv)} \right\|_{\tilde{\Phi}(\varepsilon v)} \rek w_{0}(x) dx \rge \leq {\Phi}^{-1}(1/\varepsilon)   	
    \een
    holds.
    \item  There exists a constant $B>0,$ such that for all $r>0, \varepsilon>0,$ the estimate
    \ben
    \Psi^{-1}\lge \int_{0}^{r}\Psi(\varepsilon)w_0(x)dx\rge \leq \Phi^{-1}\lge \Phi(B\varepsilon)(Iv)(r)\rge
    \een   
		holds.
  	\end{enumerate} }

	\section{Norm Inequalities}
	
  We begin this section by giving sufficient conditions so that 
norm inequalities hold for Dunkl-Hausdorff operator in Orlicz spaces.	
   
   \begin{theorem} \label{th31}
   	Suppose $\Phi$ and $\Psi$ are N- functions such that $\Phi, \tilde{\Phi} \in \Delta_{2}$ and $\tilde{\Psi} \prec \tilde{\Phi}.$ Let $v,w_0$ be weight functions with $(Iv)(\infty) = \infty$ and for $\alpha \in \mathbb{R},$  the condition \eqref{a21} holds for some constants $C_{i}, ~i = 1,2,3,4.$  If there are constant $B_1, B_2, B_3$ such that for all $r>0, \varepsilon>0$ the following hold
   	 \begin{align} \label{a41}
   		{\tilde{\Phi}}^{-1}\lge \int_{r}^{\infty} \tilde{\Phi} \lek \frac{1}{B_{1}(Iv)(x)} \left\| \frac{\chi_{(0,r)}}{\varepsilon} \right\|_{\Psi(\varepsilon w_{0})} \rek v(x) dx \rge \leq  {\tilde{\Psi}}^{-1}(1/\varepsilon)
   	\end{align}
   	\begin{align} \label{a42}
   		{\tilde{\Phi}}^{-1}\lge \int_{0}^{r} \tilde{\Phi} \lek \frac{x}{B_{2}(Iv)(x)} \left\| \frac{\chi_{(r,\infty)}}{\varepsilon x} \right\|_{\Psi(\varepsilon w_{0})} \rek v(x) dx \rge \leq  {\tilde{\Psi}}^{-1}(1/\varepsilon)
   	\end{align}
   	\begin{align} \label{a43}
   	  {\tilde{\Phi}}^{-1}\lge \int_{r}^{\infty} \tilde{\Phi} \lek \frac{1}{B_{3}(Iv)(x)} \left\| \frac{\ln(r,.)\chi_{(0,r)}} {\varepsilon} \right\|_{\Psi(\varepsilon w_{0})} \rek v(x) dx \rge \leq {\tilde{\Psi}}^{-1}(1/\varepsilon)	 
   	\end{align}
   	\begin{align} \label{a44}
   	{\tilde{\Phi}}^{-1}\lge \int_{r}^{\infty} \tilde{\Phi} \lek \frac{1}{B_{3}(Iv)(x)} \left\| \frac{\chi_{(0,r)}} {\varepsilon} \right\|_{\Psi(\varepsilon w_{0})} \ln\lr\frac{x}{r}\rr \rek v(x) dx \rge \leq {\tilde{\Psi}}^{-1}(1/\varepsilon),
   	\end{align}   	
 then for $0\leq f \downarrow$ the inequality 
   	\begin{eqnarray} \label{a45}
   		\|\mathcal{H}_{\alpha,\phi}f\|_{\Psi(w_{0})} \leq C \|f\|_{\Phi(v)}
   	\end{eqnarray}
		holds.
   \end{theorem}   
   \proof By Remark \ref{rm1}, we have that the inequality \eqref{a45}  is equivalent to 
   \be \label{a46}
   \left\|\frac{I(\mathcal{H}_{\alpha,\phi}^{*}g)}{Iv}\right\|_{\tilde{\Phi}(v)} \leq C \left\| \frac{g}{w_0} \right\|_{\tilde{\Psi}(w_{0})},   \hspace{2mm} g\geq 0.
   \ee    
   Following the method used in Theorem \ref{t1}, we see that 
   \begin{align} \label{eqn37}
   	\mathcal{I}(\mathcal{H}_{\alpha,\phi}^{*}g)(x) &\approx \int_{0}^{x} \ln\lr\frac{x}{t}\rr g(t) dt + \int_{0}^{x} g(t) dt + x\int_{x}^{\infty} \frac{g(t)}{t} dt \\ \nonumber
   	&= (\mathcal{K}g)(x) + (\mathcal{H}g)(x) + (\mathcal{S}^{*}g)(x), ~x \in [0,\infty),
   \end{align}
	where $(Sg)(x):=\frac{1}{x}\int_0^x tg(t)dt.$ Since Luxemburg norm is a Banach function norm, using the equivalence (\ref{eqn37}), the inequality \eqref{a46} is equivalent to the following three inequalities:
    \be \label{a47}
    \left\|\frac{\mathcal{A}^{*}g}{Iv}\right\|_{\tilde{\Phi}(v)} \leq C \left\| \frac{g}{w_0} \right\|_{\tilde{\Psi}(w_{0})}, 
    \ee    
    \be \label{a48}
    \left\|\frac{\mathcal{H}g}{Iv}\right\|_{\tilde{\Phi}(v)} \leq C \left\| \frac{g}{w_0} \right\|_{\tilde{\Psi}(w_{0})} 
    \ee   
		and
    \be \label{a49}
    \left\|\frac{\mathcal{K}g}{Iv}\right\|_{\tilde{\Phi}(v)} \leq C \left\| \frac{g}{w_0} \right\|_{\tilde{\Psi}(w_{0})}. 
    \ee
   Given that \eqref{a41}, \eqref{a42}, \eqref{a43} and \eqref{a44} hold, in view of Theorems C and D the following corresponding modular inequalities hold: 
    \[
   	\tilde{\Phi}^{-1}\lge \int_{0}^{\infty} \tilde{\Phi} \lek \frac{x}{(Iv)(x)} \lr  \int_{x}^{\infty} \frac{g(t)}{t} \ dt  \rr  \rek v(x) dx \rge \leq {\tilde{\Psi}}^{-1}\lge \int_{0}^{\infty} \tilde{\Psi} \lek  A \frac{g(x)}{w_{0}(x)} \rek w_{0}(x) dx \rge, 	    
    \]
     
    \[
   	\tilde{\Phi}^{-1}\lge \int_{0}^{\infty} \tilde{\Phi} \lek \frac{1}{(Iv)(x)} \lr  \int_{0}^{x} g(t) \ dt  \rr  \rek v(x) dx \rge 
   	\leq {\tilde{\Psi}}^{-1}\lge \int_{0}^{\infty} \tilde{\Psi} \lek  A\frac{g(x)}{w_{0}(x)} \rek w_{0}(x) dx \rge	
    \]
 and  
   \[
   	\tilde{\Phi}^{-1}\lge \int_{0}^{\infty} \tilde{\Phi} \lek \frac{1}{(Iv)(x)} \lr  \int_{0}^{x} \ln\lr\frac{x}{t}\rr g(t) \ dt  \rr  \rek v(x) dx \rge 
   	\leq {\tilde{\Psi}}^{-1}\lge \int_{0}^{\infty} \tilde{\Psi} \lek A \frac{g(x)}{w_{0}(x)} \rek w_{0}(x) dx \rge. 	
   \]
    Since modular inequalities imply corresponding norm inequalities in Orlicz space setting (see \cite{hm}), the inequalities \eqref{a47}, \eqref{a48} and \eqref{a49} hold and hence \eqref{a46} holds. \hfill $\square$\\
   
	\vspace{3pt}
	
    \noindent Below we give another sufficient condition for norm inequalities to hold in Orlicz space setting. 
    
   \begin{theorem} \label{th32}
   Suppose $\Phi$ and $\Psi$ are N- functions such that $\Phi, \tilde{\Phi} \in \Delta_{2}$ and $\Phi \prec \Psi.$ Let $(Iv)(\infty) = \infty$ and for $\alpha \in \mathbb{R},$  the condition \eqref{a21} holds for some constants $C_{i}, ~i = 1,2,3,4.$  If there are constants $A_1, A_2, A_3$ such that for all $r>0, \varepsilon>0$ the expressions \eqref{a23},  \eqref{a24}, \eqref{a25} and \eqref{a26} hold with $w_1=1,$ then for $0\leq f \downarrow$ the inequality
   \eqref{a45} holds.
   \end{theorem}

    \proof We know that \eqref{a45} is equivalent to \eqref{a46}. 
    Further, on using the fact (see \cite{dhk}) that 
    \ben
    \| \gamma T g\|_{\tilde{\Phi}(v)} \leq C \| \beta g\|_{\tilde{\Psi}(w)} \Longleftrightarrow \left\| \frac{T^{*}h} {w\beta}\right\|_{{\Psi}(w)} \leq C \left\| \frac{h}{\gamma v}\right\|_{{\Phi}(v)},
    \een
    where$\gamma(x)$ and $\beta(x)$ are non-negative functions, we get that \eqref{a46} is equivalent to the inequality
    \be \label{a413}
    \left\|\mathcal{H}_{\alpha,\phi}(I^{*}h)\right\|_{\Psi(w_{0})} \leq C \left\| \frac{(Iv)h}{v} \right\|_{\Phi(v)},   \hspace{2mm} h \geq 0.
    \ee    
    Thus, to prove \eqref{a45}, it is enough to prove \eqref{a413}. Using \eqref{a28}, \eqref{a211} and \eqref{a214}, we see that
    \eqref{a413} is equivalent to the following three norm inequalities:
     \be \label{a414}
    \left\|\mathcal{S}h\right\|_{\Psi(w_{0})} \leq C \left\| \frac{(Iv)h}{v} \right\|_{\Phi(v)}, 
    \ee    
    \be \label{a415}
    \left\|\mathcal{H}^{*}h\right\|_{\Psi(w_{0})} \leq C \left\| \frac{(Iv)h}{v} \right\|_{\Phi(v)}  
    \ee    
		and
    \be \label{a416}
    \left\|\mathcal{K}^{*}h\right\|_{\Psi(w_{0})} \leq C \left\| \frac{(Iv)h}{v} \right\|_{\Phi(v)},
    \ee
    where $(\mathcal{S}h)(x) = \frac{1}{x} \int_{0}^{x} t h(t) dt, ~(\mathcal{H}^{*}h)(x) = \int_{x}^{\infty} h(t) dt $ and $ (\mathcal{K}^{*}h)(x) = \int_{x}^{\infty} \ln\lr\frac{t}{x}\rr h(t) dt.$ \\
    Since \eqref{a23}, \eqref{a24}, \eqref{a25} and \eqref{a26} hold with $w_1=1,$ the corresponding modular  inequalities \eqref{a215}, \eqref{a216} and \eqref{a217} hold, and hence \eqref{a414}, \eqref{a415} and \eqref{a416} hold, which leads to that \eqref{a413} holds. \hfill $\square$ 
     
   \vspace{3pt} 
	
		In the reverse direction, now we give a situation when we may start with a norm inequality leading to modular inequality. Precisely, analogous to Proposition E, below we give an equivalence between weighted modular inequality and a certain Orlicz-Luxemburg norm inequality.
 
    \begin{theorem}
     Suppose $\Phi$ and $\Psi$ are N- functions such that $\Phi, \tilde{\Phi} \in \Delta_{2}$ and $\Phi \prec \Psi.$ Let  $w_0, v$ be weights with $(Iv)(\infty) = \infty$ and for $\alpha \in \mathbb{R},$ the condition \eqref{a21} holds for some constants $C_{i}, ~i = 1,2,3,4.$ Then for every $\varepsilon>0,$ the inequality
     \be \label{a417}
     \left\|\mathcal{H}_{\alpha,\phi}f\right\|_{\Psi(\varepsilon_{\Psi}w_0)} \leq C \|f\|_{\Phi(\varepsilon_{\Phi}v)}
     \ee
     holds for $0\leq f \downarrow$  if and only if there are constants $A_1, A_2, A_3$ such that \eqref{a23}, \eqref{a24}, \eqref{a25} and \eqref{a26} hold with $w_1=1,$ where 
     \[\varepsilon_{\Psi}=1/\Psi(1/\varepsilon) , \hspace{2mm} \varepsilon_{\Phi}=1/\Phi(1/\varepsilon).  \]
    \end{theorem}
 
  \proof Taking $w_0(x) dx = d\nu(x)$ and $v(x) dx = d\mu(x)$ and using Proposition E, the inequality \eqref{a417} is equivalent to 
   \be \label{a418}
  	{\Psi}^{-1}\lge \int_{0}^{\infty} \Psi \lek (\mathcal{H}_{\alpha,\phi} f)(x) \rek w_0(x) dx \rge \leq {\Phi}^{-1}\lge \int_{0}^{\infty} \Phi \lek  Cf(x) \rek v(x) dx \rge 
   \ee
    for $0\leq f \downarrow.$ Now, using Theorem \ref{t1} with $w_1=1,$ \eqref{a418} holds if and only if \eqref{a23}, \eqref{a24}, \eqref{a25} and \eqref{a26} hold. Hence the result follows. \hfill $\square$ 
	
		\vspace{5pt}
		
	\noindent {\bf Acknowledgment.} The first author acknowledges  the research fellowship award Ref No.: 1004 (CSIR-UGC NET June 2018) of University Grants Commission (UGC), India. The corresponding author acknowledges the Research Grant No. 02011/14/2023 NBHM(RP)/R\&D II/5951 of National Board for Higher Mathematics (NBHM), India.
		
		\vspace{3pt}
		
	\noindent {\it Conflict of interest statement.} On behalf of all authors, the corresponding author states that there is no conflict of interest.

%\bigskip
%\newpage

\vspace{15pt}

\noindent  Department of Mathematics\\
University of Delhi, Delhi - 110007, India\\
(Email: megha.madan110@gmail.com)

\bigskip
\noindent Department of Mathematics\\
Dyal Singh College (University of Delhi)\\
Lodhi Road, Delhi - 110003, India\\
(Email : arunpalsingh@dsc.du.ac.in)

\bigskip

\noindent Department of Mathematics\\
Lady Shri Ram College For Women (University of Delhi)\\
Lajpat Nagar, Delhi - 110024, India\\
(Email : monikasingh@lsr.du.ac.in)

\end{document}